\documentclass{amsart}
\usepackage{amsfonts}
\usepackage{amsmath}
\usepackage{amscd}

\input tcilatex
\theoremstyle{definition}
\theoremstyle{remark}
\numberwithin{equation}{section}

\begin{document}
\title{Positive Definite Germs of Quantum Stochastic Processes}
\author{Viacheslav Belavkin}
\address{Mathematics Department, University of Nottingham,\\
NG7 2RD, UK.}
\date{2 October 1995}
\thanks{ Published in: \textit{C. R. Acad. Sc. Paris}, \textbf{322} (1)
385--390 (1996)}
\subjclass{}
\keywords{quantum stochastic germs, completely positive flows,
conditionall-positive generators }
\maketitle

\begin{abstract}
A new notion of stochastic germs for quantum processes is introduced and a
characterisation of the stochastic differentials for positive definite (PD)
processes is found in terms of their germs for arbitrary It\^{o} algebra. A
representation theorem, giving the pseudo-Hilbert dilation for the
germ-matrix of the differential, is proved. This suggests the general form
of quantum stochastic evolution equations with respect to the Poisson
(jumps), Wiener (diffusion) or general quantum noise.
\end{abstract}

\section*{\textbf{Germes positivement definis de processus quantiques
stochastiques}}

\setcounter{equation}{1}

\begin{quote}
{\footnotesize \textbf{R\'{e}sum\'{e}.} On trouve une characterisation des
diff\'{e}rentielles stochastiques des processus quantiques positifs d\'{e}%
finis (PD) pour une alg\`{e}bre de It\^{o} arbitraire. On d\'{e}montre un th%
\'{e}or\`{e}me de repr\'{e}sentation qui donne la dilatation
pseudo-Hilbertienne de la matrice des germes de la diff\'{e}rentielle. Ceci
sugg\`{e}re une forme g\'{e}n\'{e}rale des \'{e}quations d'\'{e}volution
quantiques stochastiques par rapport aux sauts de Poisson, \`{a} la
diffusion de Wiener ou aux bruits quantiques g\'{e}n\'{e}ralis\'{e}s.}
\end{quote}

\subsection*{\emph{Version fran\c{c}aise abr\'{e}g\'{e}e}.}

\subsubsection*{1. Introduction.}

Le but du pr\'{e}sent article est de g\'{e}n\'{e}raliser le th\'{e}or\`{e}me
de Evans-Lewis \cite{EvL} sur la construction diff\'{e}rentielle de la
dilatation de Stinespring \cite{Stn} au cas des diff\'{e}rentielles
stochastiques parametris\'{e}es par une $\star $-alg\`{e}bre quelconque.
Comme j'ai d\'{e}montr\'{e} dans \cite{Bcs}, tout processus stochastique
stationaire $\Lambda \left( a\right) $ \`{a} incr\'{e}ments ind\'{e}pendents
peut \^{e}tre repr\'{e}sent\'{e} dans un espace de Fock $\mathcal{F}$ sur
l'espace $\mathcal{K}\otimes L^2\left( \mathbb{R}_{+}\right) $ des fonctions
au carr\'{e} int\'{e}grable sur $\mathbb{R}_{+}$ \`{a} valeurs dans $%
\mathcal{K}$ sous la forme $\mathrm{A}\left( \mathbf{a}\right) =a_\nu ^\mu 
\mathrm{A}_\mu ^\nu $, o\`{u} 
\begin{equation*}
a_\nu ^\mu \mathrm{A}_\mu ^\nu \left( t\right) =a_{\bullet }^{\bullet }%
\mathrm{A}_{\bullet }^{\bullet }\left( t\right) +a_{+}^{\bullet }\mathrm{A}%
_{\bullet }^{+}\left( t\right) +a_{\bullet }^{-}\mathrm{A}_{-}^{\bullet
}\left( t\right) +a_{+}^{-}\mathrm{A}_{-}^{+}\left( t\right) ,
\end{equation*}
est la d\'{e}composition canonique de $\Lambda $ dans l'\'{e}change des
processus $\mathrm{A}_{\bullet }^{\bullet }$ de cr\'{e}ation $\mathrm{A}%
_{\bullet }^{+}$, d'annihilation $\mathrm{A}_{-}^{\bullet }$ et temps $%
\mathrm{A}_{-}^{+}=t\mathrm{I}$ par rapport \`{a} l'\'{e}tat vide dans $%
\mathcal{F}$. Ainsi l'alg\`{e}bre de parametrisation $\mathfrak{a}$ peut 
\^{e}tre identifi\'{e}e \`{a} l'alg\`{e}bre des quadruplets GNS $\boldsymbol{%
a}=\left( a_\nu ^\mu \right) _{\nu =+,\bullet }^{\mu =-,\bullet }$, associ%
\'{e}s \`{a} la $*$ -functionelle lin\'{e}aire positive $l\left( a\right)
=\langle \Lambda \left( t,a\right) \rangle /t$.

Les principaux r\'{e}sultats de cet article sont des th\'{e}or\`{e}mes de
characterisation et de repr\'{e}sentation des germes stochastiques
quantiques $\boldsymbol{\lambda }=\boldsymbol{i}+\boldsymbol{\alpha }$ des $%
* $-repr\'{e}sentations unitaires $i$ qui d\'{e}finissent les diff\'{e}%
rentielles de It\^{o} $\mathrm{d}\phi =\alpha _\nu ^\mu \mathrm{dA}_\mu ^\nu 
$ \`{a} $t=0$ pour les processus compl\`{e}tement positifs $\phi _t $ avec $%
\phi _0=i$.

\subsubsection*{2. Notations.}

Pour tout espace pr\'{e}-Hilbertien $\mathcal{D}$, le symbole $\mathcal{D}%
^{\prime }$ d\'{e}signe l'espace dual de toutes les fonctionelles antilin%
\'{e}aires $\eta ^{\prime }:\eta \mapsto \langle \eta |\eta ^{\prime
}\rangle $ avec le plongement naturel $\mathcal{D}\subseteq \mathcal{D}%
^{\prime },$ donn\'{e} par le produit int\'{e}rieur $\langle \eta ^{\prime
}|\eta \rangle =\langle \eta |\eta ^{\prime }\rangle ^{*}$ pour $\eta
^{\prime }\in \mathcal{D}$. Toute forme sesquilin\'{e}aire $\langle \eta
|Q\eta \rangle $ sur $\mathcal{D}$ s'identifie \`{a} l'op\'{e}rateur lin\'{e}%
aire $Q:\mathcal{D}\rightarrow \mathcal{D}^{\prime }$. L'espace lin\'{e}aire 
$\mathcal{Q}\left( \mathcal{D}\right) $ de tous ces op\'{e}rateurs est muni
d'une involution $Q\mapsto Q^{\dagger },$ o\`{u} $Q^{\dagger }:\mathcal{D}%
\rightarrow \mathcal{D}^{\prime }$ est l'op\'{e}rateur Hermitien adjoint: $%
\langle \eta |Q^{\dagger }\eta \rangle =\langle \eta |Q\eta \rangle ^{*}.$
Tout op\'{e}rateur lin\'{e}aire $A:\mathcal{D\rightarrow D}$ avec son
adjoint $A^{\dagger }:\eta \mapsto A^{\dagger }\eta \in \mathcal{D}$ pour
tout $\eta \in \mathcal{D}$ peut \^{e}tre continu\'{e} \`{a} $\mathcal{D}%
^{\prime }$ comme $A^{\prime }=A^{\dagger *},$ o\`{u} $A^{*}:\mathcal{D}%
^{\prime }\rightarrow \mathcal{D}^{\prime }$ est l'op\'{e}rateur dual defini
par $\langle \eta |A^{*}\eta ^{\prime }\rangle =\langle A\eta |\eta ^{\prime
}\rangle $. La $\star $-alg\`{e}bre de tous ces op\'{e}rateurs dans $%
\mathcal{D}$ est not\'{e}e $\mathcal{A}\left( \mathcal{D}\right) .$

Soient $\mathfrak{D}$ la $\mathcal{D}$-envelope des vecteurs coh\'{e}rents
(exponentiels) de l'espace de Fock $\mathcal{F}$, et $\left( \mathfrak{D}%
_{t}\right) _{t>0}$ la filtration naturelle des sous-espaces $\mathfrak{D}%
_{t}$ engendr\'{e} par l'espace pr\'{e}-Hilbertien $\mathcal{D}=\mathfrak{D}%
_{0}$ et les vecteurs coh\'{e}rents $f^{\otimes }$ sur les fonctions au carr%
\'{e} integrable $f^{\bullet }:[0,t)\rightarrow \mathcal{K}$, o\`{u} $%
\mathcal{K}$ est l'espace pr\'{e}-Hilbertien associ\'{e} \`{a} l'alg\`{e}bre
de It\^{o} $\mathfrak{a}$. Soit $B$ un semigroupe \`{a} involution $x\mapsto
x^{\star },\left( x^{\star }y\right) ^{\star }=y^{\star }x$. Le processus PD
quantique stochastique adapt\'{e} sur un $\ast $-semigroup $B$ \`{a} unit%
\'{e} est d\'{e}crit par une famille de param\'{e}tres $\phi =\left( \phi
_{t}\right) _{t>0}$ d'applications positives d\'{e}finies $\phi _{t}:$ $%
B\rightarrow \mathcal{Q}\left( \mathfrak{D}\right) $, 
\begin{equation}
\sum_{x,y\in B}\langle \mathfrak{h}_{y}|\phi _{t}\left( y^{\star }x\right) 
\mathfrak{h}_{x}\rangle :=\sum_{k,l}\langle \mathfrak{h}^{l}|\phi _{t}\left(
x_{l}^{\star }x_{k}\right) \mathfrak{h}^{k}\rangle \geq 0  \tag{0.1}
\label{2.1}
\end{equation}
pour toute suite $\mathfrak{h}^{k}\in \mathfrak{D}$ et $x_{k}\in B$, $%
k=1,2,...$, satisfaisant la propri\'{e}t\'{e} de causalit\'{e} 
\begin{equation*}
\phi _{t}\left( x\right) \left( \eta \otimes f^{\otimes }\right) =\eta
^{\prime }\otimes f^{\otimes }\quad \forall x\in B,\eta \in \mathfrak{D}%
_{t},f^{\bullet }\in \mathcal{K}\otimes L^{2}[t,\infty ),
\end{equation*}
o\`{u} $\eta ^{\prime }\in \mathfrak{D}_{t}^{\prime }$. La famille
stochastistiquement diff\'{e}rentiable $\phi $ par rapport \`{a} un
processus quantique stationaire $\Lambda $ \`{a} increments ind\'{e}pendents
engendr\'{e} par l'alg\`{e}bre de It\^{o} $\mathfrak{a}$ est donn\'{e} par
l'int\'{e}grale quantique stochastique \cite{HdP},\cite{Bcs}, 
\begin{equation}
\phi _{t}\left( x\right) -\phi _{0}\left( x\right) =\int_{0}^{t}\boldsymbol{%
\alpha }\left( r,x\right) \mathrm{d}\mathbf{A}:=\sum_{\mu ,\nu
}\int_{0}^{t}\alpha _{\nu }^{\mu }\left( r,x\right) \mathrm{d}A_{\mu }^{\nu
},\qquad \text{ }x\in B\qquad  \tag{0.2}  \label{2.2}
\end{equation}
o\`{u} $\boldsymbol{\alpha }\left( t,x\right) $ pour tout $x\in B$ est un
processus quantique stochastique adapt\'{e} \`{a} valeurs $\boldsymbol{a}%
\left( t\right) =\langle \mathfrak{h}|\boldsymbol{\alpha }\left( t,x\right) 
\mathfrak{h}\rangle $ dans les quadruplets $\left( a_{\nu }^{\mu }\right)
\left( t\right) $ qui repr\'{e}sentent l'alg\`{e}bre $\mathfrak{a}$. Les int%
\'{e}grateurs quantiques stochastiques $\mathbf{A}\left( t\right) =\left( 
\mathrm{A}_{\mu }^{\nu }\right) _{\mu =-,\bullet }^{\nu =+,\bullet }\left(
t\right) $ sont formellement d\'{e}finis par la $\flat $-propri\'{e}t\'{e} 
\cite{Bcs} $\left( \boldsymbol{a}\mathbf{A}\right) ^{\dagger }=\boldsymbol{a}%
^{\flat }\mathbf{A}$ et la table de multiplication de Hudson-Parthasarathy 
\cite{HdP} 
\begin{equation}
\left( \boldsymbol{a}\mathrm{d}\mathbf{A}\right) \left( \boldsymbol{b}%
\mathrm{d}\mathbf{A}\right) =\left( \boldsymbol{a}\boldsymbol{b}\right) 
\mathrm{d}\mathbf{A},\qquad  \tag{0.3}  \label{2.3}
\end{equation}
o\`{u} $a_{-\nu }^{\flat \mu }=a_{-\mu }^{\nu \dagger }$ par rapport \`{a}
la sym\'{e}trie $-(-)=+,$ $-(+)=-$ des indices $\left( -,+\right) $
seulement, et $\left( \boldsymbol{b}\boldsymbol{a}\right) _{\nu }^{\mu
}=b_{\bullet }^{\mu }a_{\nu }^{\bullet }$.

\subsubsection*{3. R\'{e}sultats}

\begin{theorem}
Supposons que le PD processus $\phi _{t}$ poss\`{e}de la diff\'{e}rentielle
quantique stochastique $\boldsymbol{\alpha }\left( t\right) \mathrm{d}%
\mathbf{A}$ \`{a} $t\geq 0$. Alors la matrice des germes $\boldsymbol{%
\lambda }$=$\boldsymbol{\phi }_{t}+\boldsymbol{\alpha }\left( t\right) $
pour $\boldsymbol{\phi }_{t}\left( x\right) =\left( \phi _{t}\left( x\right)
\delta _{\nu }^{\mu }\right) _{\nu =+,\bullet }^{\mu =-,\bullet }$ est
conditionnellement positive d\'{e}finite 
\begin{equation*}
\sum_{k,l}\langle \boldsymbol{\eta }^{l}|\boldsymbol{\iota }\left(
x_{l}^{\star }x_{k}\right) \boldsymbol{\eta }^{k}\rangle =0\Rightarrow
\sum_{k,l}\langle \boldsymbol{\eta }^{l}|\boldsymbol{\lambda }\left(
x_{l}^{\star }x_{k}\right) \boldsymbol{\eta }^{k}\rangle \geq 0,
\end{equation*}
o\`{u} $\boldsymbol{\eta }^{k}\in \mathfrak{D}_{t}\oplus \mathfrak{D}%
_{t}^{\bullet }$, $\mathfrak{\ D}_{t}^{\bullet }=\mathcal{K}\otimes 
\mathfrak{D}_{t}$, par rapport \`{a} l'application d\'{e}g\'{e}ner\'{e}e $%
\boldsymbol{\iota }=\left( \iota _{\nu }^{\mu }\right) _{\nu \neq -}^{\mu
\neq +},$ $\iota _{\nu }^{\mu }\left( x\right) =\phi _{t}\left( x\right)
\delta _{\nu }^{+}\delta _{-}^{\mu }$, toutes les deux \'{e}crites sous la
forme matricielle suivante 
\begin{equation}
\boldsymbol{\lambda }=\left( 
\begin{array}{cc}
\lambda & \lambda _{\bullet } \\ 
\lambda ^{\bullet } & \lambda _{\bullet }^{\bullet }%
\end{array}
\right) ,\qquad \text{ }\boldsymbol{\iota }=\left( 
\begin{array}{cc}
\iota & 0 \\ 
0 & 0%
\end{array}
\right) \qquad  \tag{0.4}  \label{2.5}
\end{equation}
avec $\lambda =\alpha _{+}^{-}\left( t\right) ,\quad $ $\lambda ^{\bullet
}=\alpha _{+}^{\bullet }\left( t\right) ,\quad $ $\lambda _{\bullet }=\alpha
_{\bullet }^{-}\left( t\right) ,\quad \lambda _{\bullet }^{\bullet }=\phi
_{t}\delta _{\bullet }^{\bullet }+\alpha _{\bullet }^{\bullet }\left(
t\right) $, et $\iota =\phi _{t}$.
\end{theorem}

\begin{theorem}
Les conditions suivantes sont \'{e}quivalentes:

\begin{enumerate}
\item[(i)] L'application de germes $\boldsymbol{\lambda }=\boldsymbol{i}+%
\boldsymbol{\alpha }$ \`{a} $t=0$ est conditionellemnt positive d\'{e}finie
par rapport \`{a} la repr\'{e}sentation matricielle $\boldsymbol{\iota }$ in
(\ref{2.5}), o\`{u} $\iota =i$.

\item[(ii)] Il \'{e}xiste un espace pr\'{e}-Hilbertien $\mathcal{D}^{\circ } 
$, une $\ast $- repr\'{e}sentation \`{a} unit\'{e} $j$ de $B$ dans $\mathcal{%
A}\left( \mathcal{D}^{\circ }\right) $, 
\begin{equation}
j\left( y^{\star }x\right) =j\left( y\right) ^{\ast }j\left( x\right) ,\quad
j\left( 1\right) =I,  \tag{0.6}  \label{3.6}
\end{equation}
une $\left( j,i\right) $-d\'{e}rivation de $B$, 
\begin{equation}
k\left( y^{\star }x\right) =j\left( y\right) ^{\ast }k\left( x\right)
+k\left( y^{\star }\right) i\left( x\right) ,  \tag{0.7}  \label{3.7}
\end{equation}
\`{a} valeurs dans les op\'{e}rateurs $\mathcal{D}\rightarrow \mathcal{D}%
^{\circ }$, et une application $l:B\rightarrow \mathcal{Q}\left( \mathcal{D}%
\right) $ telle que 
\begin{equation}
l\left( y^{\star }x\right) =i\left( y\right) ^{\ast }l\left( x\right)
+l\left( y^{\star }\right) i\left( x\right) +k\left( y\right) ^{\dagger
}k\left( x\right) ,  \tag{0.8}  \label{3.8}
\end{equation}
avec l'application adjointe $l\left( x\right) ^{\dagger }=l\left( x^{\star
}\right) +Di\left( x^{\star }\right) -i\left( x\right) ^{\ast }D\quad $
telle que $l\left( x\right) +Di\left( x\right) =\lambda \left( x\right)
=l^{\ast }\left( x\right) +i\left( x\right) ^{\prime }D$, 
\begin{equation*}
L_{\circ }^{\bullet }k\left( x\right) +L_{+}^{\bullet }i\left( x\right)
=\lambda ^{\bullet }\left( x\right) ,\quad \lambda _{\bullet }\left(
x\right) =k^{\ast }\left( x\right) L_{\bullet }^{\circ }+i\left( x\right)
^{\prime }L_{\bullet }^{-},
\end{equation*}
et $\lambda _{\bullet }^{\bullet }\left( x\right) =L_{\circ }^{\bullet
}j\left( x\right) L_{\bullet }^{\circ }$ pour certains op\'{e}rateurs $%
L_{\circ }^{\bullet }:\mathcal{D}^{\circ }\rightarrow \mathcal{D}^{\bullet }$
avec les adjoints $L_{\bullet }^{\circ }=L_{\circ }^{\bullet \dagger }$ and $%
L_{+}^{\bullet }:\mathcal{D}\rightarrow \mathcal{D}^{\bullet }$ o\`{u} $%
L_{\bullet }^{-}=L_{+}^{\bullet \dagger }$
\end{enumerate}
\end{theorem}

\setcounter{section}{0}

\section{Introduction.}

The purpose of this paper is to extend the Evans-Lewis theorem \cite{EvL}
for the differential construction of the Stinespring dilation \cite{Stn} to
the stochastic differentials, generated by an It\^{o} $\star $--algebra 
\begin{equation*}
\mathrm{d}\Lambda \left( a\right) ^{\dagger }\mathrm{d}\Lambda \left(
a\right) =\mathrm{d}\Lambda \left( a^{\star }a\right) ,\quad \sum \lambda _i%
\mathrm{d}\Lambda \left( a_i\right) =\mathrm{d}\Lambda \left( \sum \lambda
_ia_i\right) ,\quad \mathrm{d}\Lambda \left( a\right) ^{\dagger }=\mathrm{d}%
\Lambda \left( a^{\star }\right) ,
\end{equation*}
of independent increments $\mathrm{d}\Lambda \left( t,a\right) =\Lambda
\left( t+\mathrm{d}t,a\right) -\Lambda \left( t,a\right) $, with given mean
values $\langle \mathrm{d}\Lambda \left( t,a\right) \rangle =l\left(
a\right) \mathrm{d}t$, $a\in \mathfrak{a}$. Here $l:\mathfrak{a}\rightarrow 
\mathbb{C}$ is a positive linear functional on the parametrizing $\star $%
-algebra $\mathfrak{a}$, defining the GNS representation $a\mapsto 
\boldsymbol{a}=$ $\left( a_\nu ^\mu \right) _{\nu =+,\bullet }^{\mu
=-,\bullet }$ of $\mathfrak{a}$ in terms of the quadruples 
\begin{equation}
a_{\bullet }^{\bullet }=j\left( a\right) ,\quad a_{+}^{\bullet }=k\left(
a\right) ,\quad a_{\bullet }^{-}=k^{*}\left( a\right) ,\quad
a_{+}^{-}=l\left( a\right) ,
\end{equation}
where $j\left( a^{\star }a\right) =j\left( a\right) ^{*}j\left( a\right) $
is the operator representation $j\left( a\right) ^{*}k\left( a\right)
=k\left( a^{\star }a\right) $ on the pre-Hilbert space $\mathcal{K}$ of the
Kolmogorov decomposition $l\left( a^{*}a\right) =k\left( a\right)
^{*}k\left( a\right) $, and $k^{*}\left( a\right) =k\left( a^{\star }\right)
^{*}.$

As was proved in \cite{Bcs}, the stochastic process $\Lambda \left( a\right) 
$ with independent increments can be represented in the Fock space $%
\mathfrak{F}$ over the space of $\mathcal{K}$ -valued square-integrable
functions on $\mathbb{R}_{+}$ as $\mathrm{A}\left( \mathbf{a}\right) =a_\nu
^\mu \mathrm{A}_\mu ^\nu $, where 
\begin{equation}
a_\nu ^\mu \mathrm{A}_\mu ^\nu \left( t\right) =a_{\bullet }^{\bullet }%
\mathrm{A}_{\bullet }^{\bullet }\left( t\right) +a_{+}^{\bullet }\mathrm{A}%
_{\bullet }^{+}\left( t\right) +a_{\bullet }^{-}\mathrm{A}_{-}^{\bullet
}\left( t\right) +a_{+}^{-}\mathrm{A}_{-}^{+}\left( t\right) ,
\end{equation}
is the canonical decomposition of $\Lambda $ into the exchange $\mathrm{A}%
_{\bullet }^{\bullet }$, creation $\mathrm{A}_{\bullet }^{+}$, annihilation $%
\mathrm{A}_{-}^{\bullet }$ and time $\mathrm{A}_{-}^{+}=t\mathrm{I}$
processes with respect to the vacuum state in $\mathfrak{F}$. Thus the
parametrizing algebra $\mathfrak{a}$ can be identified with the $\flat $%
-algebra of quadruples $\boldsymbol{a}$ with respect to the product $\left( 
\boldsymbol{b}\boldsymbol{a}\right) _\nu ^\mu =b_{\bullet }^\mu a_\nu
^{\bullet }$ and the involution $a_{-\nu }^{\flat \mu }=a_{-\mu }^{\nu
\dagger }$, where $-(-)=+$, $-\bullet =\bullet $, $-(+)=-$.

The main results of this paper are characterization and representation
theorems for the quantum stochastic germs $\boldsymbol{\lambda }=\boldsymbol{%
i}+\boldsymbol{\alpha }$ of unital $*$-representations $i$ defining the It%
\^{o} differentials $\mathrm{d}\phi =\alpha _\nu ^\mu \mathrm{dA}_\mu ^\nu $
at $t=0$ for completely positive processes $\phi _t$ with $\phi _0=i$. The
Evans-Lewis case $\Lambda \left( t,a\right) =at\mathrm{I}$ is described by
the simplest one-dimensional algebra $\mathfrak{a}$ with $l\left( a\right)
=a $ and the nilpotent multiplication $a^{\star }a=0$ corresponding to the
non-stochastic (Newton) calculus $\left( \mathrm{d}t\right) ^2=0$. The
unital $\star $-algebra $\mathfrak{a}=\mathbb{C}$ with the commutative
multiplication $a^{\star }a=\left| a\right| ^2$ has the GNS representation $%
a_\nu ^\mu =a$, corresponding to $\Lambda \left( t,a\right) =a\mathrm{P}%
\left( t\right) $, where $\mathrm{P}$ is the standard Poisson process $%
\mathrm{P}=\sum \mathrm{A}_\mu ^\nu $. The standard Wiener process $\mathrm{Q%
}=\mathrm{A}_{-}^{\bullet }+\mathrm{A}_{\bullet }^{+}$ is described by the
second order nilpotent algebra $\mathfrak{a}$ of quadruples $%
a_{+}^{-}=a,\quad a_{\bullet }^{-}=b=a_{+}^{\bullet },\quad a_{\bullet
}^{\bullet }=0$, corresponding to $\Lambda \left( t,a\right) =at\mathrm{I}+b%
\mathrm{Q}\left( t\right) $. Thus our results are applicable also to the
classical stochastic differentials of completely positive processes,
corresponding to the commutative It\^{o} algebras, which are decomposable
into the Wiener, Poisson and Newton orthogonal components.

\section{Notation.}

Throughout the pre-Hilbert space $\mathcal{D}$ is complex, $\mathcal{D}%
^{\prime }$ denotes the dual space of all antilinear functionals $\eta
^{\prime }:\eta \mapsto \langle \eta |\eta ^{\prime }\rangle $ with the
natural embedding $\mathcal{D}\subseteq \mathcal{D}^{\prime }$, $%
\left\langle \eta |\eta ^{\prime }\right\rangle =\left\| \eta \right\| ^2$
if $\eta ^{\prime }=\eta \in \mathcal{D}$, $\mathcal{Q}\left( \mathcal{D}%
\right) $ denotes the space of all linear operator $Q:\mathcal{D}\rightarrow 
\mathcal{D}^{\prime }$ identified with the sesquilinear forms $\langle \eta
|Q\eta \rangle $ on $\mathcal{D}$ and $\mathcal{A}\left( \mathcal{D}\right) $
denotes the unital $\star $-algebra of all operators $Q\in \mathcal{Q}\left( 
\mathcal{D}\right) $ with $Q\mathcal{D}\subseteq \mathcal{D}$, $Q^{\dagger }%
\mathcal{D\subseteq D}$, where $\langle \eta |Q^{\dagger }\eta \rangle
=\langle \eta |Q\eta \rangle ^{*}$. Any operator $A\in \mathcal{A}\left( 
\mathcal{D}\right) $ can be extended onto $\mathcal{D}^{\prime }$ as $%
A^{\prime }=A^{\dagger *}$, where $A^{*}:\mathcal{D}^{\prime }\rightarrow 
\mathcal{D}^{\prime }$ is the dual operator, defined by $\langle \eta
|A^{*}\eta ^{\prime }\rangle =\langle A\eta |\eta ^{\prime }\rangle $.

Let $\mathfrak{D}$ denote the $\mathcal{D}$-span of the coherent
(exponential) vectors $f^{\otimes }:\tau \mapsto \otimes _{t\in \tau
}f^{\bullet }\left( t\right) $, $f^{\bullet }\in \mathfrak{K}$ of the Fock
space $\mathfrak{F}$ over $\mathfrak{K}=\mathcal{K}\otimes L^2\left( \mathbb{%
R}_{+}\right) $, where $\mathcal{K}$ is the pre-Hilbert space, associated
with the It\^{o} algebra $\mathfrak{a}$, and $\left( \mathfrak{D}_t\right)
_{t>0}$ be the natural filtration of the subspaces $\mathfrak{D}_t$
generated by $\eta \otimes f^{\otimes }$, $\eta \in \mathcal{D}$ and $%
f^{\otimes }$, $f^{\bullet }\in \mathcal{K}\otimes L^2[0,t) $, embedded into 
$\mathfrak{F}$ as $f^{\otimes }\left( \tau \right) =0$ for any finite $\tau
\subset \mathbb{R}_{+}$ if $\tau \cap [t,\infty )\neq \emptyset $. Let $B$
be a unital semigroup with involution $x\mapsto x^{\star }$, $\left(
x^{\star }y\right) ^{\star }=y^{\star }x$. Say $B=\mathfrak{b}$ is the $%
\star $-semigroup $1+\mathfrak{a}$ of the It\^{o} algebra $\mathfrak{a}$
with $\left( 1+a\right) ^{\star }\left( 1+b\right) =1+a\star b$, where $%
a\star b=b+a^{\star }b+a^{\star }$, or $B=\mathcal{B}$ is an operator
algebra $\mathcal{B\subseteq A}\left( \mathcal{D}\right) $. The quantum
stochastic adapted PD process over $B\mathcal{\ }$is described by a one
parameter family $\phi =\left( \phi _t\right) _{t>0}$ of positive definite
maps $\phi _t:$ $B\rightarrow \mathcal{Q}\left( \mathfrak{D}\right) $, 
\begin{equation}
\sum_{x,y\in B}\langle \mathfrak{h}_y|\phi _t\left( y^{\star }x\right) 
\mathfrak{h}_x\rangle :=\sum_{k,l}\langle \mathfrak{h}^l|\phi _t\left(
x_l^{\star }x_k\right) \mathfrak{h}^k\rangle \geq 0  \label{2.1}
\end{equation}
for any sequence $\mathfrak{h}^k\in \mathfrak{D}$ and $x_k\in B,k=1,2,...$,
satisfying the causality property 
\begin{equation*}
\phi _t\left( x\right) \left( \eta \otimes f^{\otimes }\right) =\eta
^{\prime }\otimes f^{\otimes }\quad \forall x\in B,\eta \in \mathfrak{D}%
_t,f^{\bullet }\in \mathcal{K}\otimes L^2[t,\infty ),
\end{equation*}
where $\eta ^{\prime }\in \mathfrak{D}_t^{\prime }$. The positive
definiteness (\ref{2.1}) obviously implies the $*$-property $\phi
_t^{*}=\phi _t$, where $\phi _t^{*}\left( x\right) =\phi _t\left( x^{\star
}\right) ^{\dagger }$. The stochastically differentiable family $\phi $ with
respect to a quantum stationary process $\Lambda $ with independent
increments generated by the It\^{o} algebra $\mathfrak{a}$ is given by the
quantum stochastic integral \cite{HdP},\cite{Bcs}, 
\begin{equation}
\phi _t\left( x\right) -\phi _0\left( x\right) =\int_0^t\boldsymbol{\alpha }%
\left( r,x\right) \mathrm{d}\mathbf{A}:=\sum_{\mu ,\nu }\int_0^t\alpha _\nu
^\mu \left( r,x\right) \mathrm{dA}_\mu ^\nu ,\qquad \text{ }x\in B\qquad
\label{2.2}
\end{equation}
where $\boldsymbol{\alpha }\left( t,x\right) $ for each $x\in B$ is a
quantum stochastic adapted process with the values $\boldsymbol{a}\left(
t\right) =\langle \mathfrak{h}|\boldsymbol{\alpha }\left( t,x\right) 
\mathfrak{h}\rangle $ into the quadruples $\left( a_\nu ^\mu \right) \left(
t\right) $, representing the algebra $\mathfrak{a}$. The quantum stochastic
integrators $\mathbf{A}\left( t\right) =\left( \mathrm{A}_\mu ^\nu \right)
_{\mu =-,\bullet }^{\nu =+,\bullet }\left( t\right) $ are formally defined
by the $\flat $-property \cite{Bcs} $\left( \boldsymbol{a}\mathbf{A}\right)
^{\dagger }=\boldsymbol{a}^{\flat }\mathbf{A}$ and the Hudson-Parthasarathy
multiplication table \cite{HdP} 
\begin{equation}
\left( \boldsymbol{a}\mathrm{d}\mathbf{A}\right) \left( \boldsymbol{b}%
\mathrm{d}\mathbf{A}\right) =\left( \boldsymbol{a}\boldsymbol{b}\right) 
\mathrm{d}\mathbf{A}.\qquad  \label{2.3}
\end{equation}
The quantum stochastic derivatives $\alpha _\nu ^\mu $ for the PD processes $%
\phi _t^{*}=\phi _t$ should obviously satisfy the $\flat $-property $%
\boldsymbol{\alpha }^{\flat }=\boldsymbol{\alpha }$ where $\alpha _{-\mu
}^{\flat \nu }=\alpha _{-\nu }^{\mu *},$ $\alpha _\nu ^{\mu *}\left(
t,x\right) =\alpha _\nu ^\mu \left( t,x^{\star }\right) ^{\dagger }$.

In order to make the formulation of the dilation theorem as concise as
possible, we need the notion of the $\flat $-representation of $B$ in the
operator algebra $\mathcal{A}\left( \mathcal{E}\right) $ of a pseudo-Hilbert
space $\mathcal{E}=\mathcal{D}^{\prime }\oplus \mathcal{D}^{\circ }\oplus 
\mathcal{D}$ with respect to the indefinite metric 
\begin{equation}
\left( \xi \right| \left. \xi \right) =2\func{Re}\left( \xi ^{-}\right|
\left. \xi ^{+}\right) +\left\| \xi ^{\circ }\right\| ^2+\left\| \xi
^{+}\right\| _D^2  \label{3.4}
\end{equation}
for the triples $\xi ^\mu ,\mu =-,\circ ,+$, where $\xi ^{+}\in \mathcal{D}$%
, $\xi ^{\circ }\in \mathcal{D}^{\circ }$, $\xi ^{-}\in \mathcal{D}^{\prime
} $, $\mathcal{D}^{\circ }$ is a pre-Hilbert space, and $\left\| \eta
\right\| _D^2=\left\langle \eta |D\eta \right\rangle $. The operators $L\in 
\mathcal{A}\left( \mathcal{E}\right) $ given by $3\times 3$-block-matrices $%
\mathbf{L=}\left[ L_\nu ^\mu \right] _{\nu =-,\circ ,+}^{\mu =-,\circ ,+}$
have the Pseudo-Hermitian adjoints $\left( \xi |L^{\flat }\xi \right)
=\left( L\xi |\xi \right) $, which are defined by the Hermitian adjoints $%
L_\nu ^{\dagger \mu }=L_\mu ^{\nu *}$ as $\mathbf{L}^{\flat }=\mathbf{G}^{-1}%
\mathbf{L}^{\dagger }\mathbf{G}$ respectively to the indefinite metric
tensor $\mathbf{G}=\left[ G_{\mu \nu }\right] $ and its inverse $\mathbf{G}%
^{-1}=\left[ G^{\mu \nu }\right] $, given by 
\begin{equation}
\mathbf{G}=\left[ 
\begin{array}{ccc}
0 & 0 & I \\ 
0 & I_{\circ }^{\circ } & 0 \\ 
I & 0 & D%
\end{array}
\right] ,\qquad \mathbf{G}^{-1}=\left[ 
\begin{array}{ccc}
-D & 0 & I \\ 
0 & I_{\circ }^{\circ } & 0 \\ 
I & 0 & 0%
\end{array}
\right]  \label{3.5}
\end{equation}
with Hermitian $D$, where $I_{\circ }^{\circ }$ is the identity operator in $%
\mathcal{D}^{\circ }$, being equal $I_{\bullet }^{\bullet }=\mathbf{1\otimes 
}I$ in the case $\mathcal{D}^{\circ }=\mathcal{D}^{\bullet }$, where $%
\mathcal{D}^{\bullet }=\mathcal{K}\otimes \mathcal{D}.$

\section{The Results}

1. The following theorem in particular proves that the quantum stochastic
germ $\boldsymbol{\lambda }=\left( \lambda _\nu ^\mu \right) _{\nu
=+,\bullet }^{\mu =-,\bullet }$, $\lambda _\nu ^\mu \left( x\right) =i\left(
x\right) \delta _\nu ^\mu +\alpha _\nu ^\mu \left( x\right) $ of a unital $*$%
-representation $\phi _0=i$ of $B$ on $\mathcal{D}$, defined by the quantum
stochastic derivatives $\boldsymbol{\alpha }=\left( \alpha _\nu ^\mu \right)
_{\nu =+,\bullet }^{\mu =-,\bullet }$ for a PD adapted process $\phi $ at $%
t=0$, must be conditionally PD with respect to the embedded representation $%
\boldsymbol{\iota }\left( x\right) =\left( i\left( x\right) \delta _\nu
^{+}\delta _{-}^\mu \right) _{\nu =+,\bullet }^{\mu =-,\bullet }$ on $%
\mathcal{D}\oplus \mathcal{D}^{\bullet }$. Here $\delta _\nu ^\mu $ is the
Kronecker delta with $\delta _{\bullet }^{\bullet }=\mathbf{1}$, where $%
\mathbf{1}$ is the identity operator on $\mathcal{K}$, and $X\delta
_{\bullet }^{\bullet }=X\otimes \mathbf{1}$.

\begin{theorem}
Suppose that the PD process $\phi _t$ has the quantum stochastic
differential $\boldsymbol{\alpha }\left( t\right) \mathrm{d}\mathbf{A}$ at a 
$t\geq 0$. Then the germ-map $\boldsymbol{\lambda }$=$\boldsymbol{\phi }_t+%
\boldsymbol{\alpha }\left( t\right) $ for $\boldsymbol{\phi }_t\left(
x\right) =\left( \phi _t\left( x\right) \delta _\nu ^\mu \right) _{\nu
=+,\bullet }^{\mu =-,\bullet }$ is conditionally positive definite 
\begin{equation*}
\sum_{k,l}\langle \boldsymbol{\eta }^l|\boldsymbol{\iota }\left( x_l^{\star
}x_k\right) \boldsymbol{\eta }^k\rangle =0\Rightarrow \sum_{k,l}\langle 
\boldsymbol{\eta }^l|\boldsymbol{\lambda }\left( x_l^{\star }x_k\right) 
\boldsymbol{\eta }^k\rangle \geq 0
\end{equation*}
where $\boldsymbol{\eta }^k\in \mathfrak{D}_t\oplus \mathfrak{D}_t^{\bullet
},\mathfrak{D}_t^{\bullet }=\mathcal{K}\otimes \mathfrak{D}_t$, with respect
to the degenerate map $\boldsymbol{\iota }=\left( \iota _\nu ^\mu \right)
_{\nu =\bullet ,+}^{\mu =-,\bullet }$, $\iota _\nu ^\mu \left( x\right)
=\phi _t\left( x\right) \delta _\nu ^{+}\delta _{-}^\mu $, both written in
the matrix form as 
\begin{equation}
\boldsymbol{\lambda }=\left( 
\begin{array}{cc}
\lambda & \lambda _{\bullet } \\ 
\lambda ^{\bullet } & \lambda _{\bullet }^{\bullet }%
\end{array}
\right) ,\qquad \text{ }\boldsymbol{\iota }=\left( 
\begin{array}{cc}
\iota & 0 \\ 
0 & 0%
\end{array}
\right) \qquad  \label{2.5}
\end{equation}
with $\lambda =\alpha _{+}^{-}\left( t\right) ,\quad $ $\lambda ^{\bullet
}=\alpha _{+}^{\bullet }\left( t\right) ,\quad $ $\lambda _{\bullet }=\alpha
_{\bullet }^{-}\left( t\right) ,\quad \lambda _{\bullet }^{\bullet }=\phi
_t\delta _{\bullet }^{\bullet }+\alpha _{\bullet }^{\bullet }\left( t\right)
,$and $\iota =\phi _t$ such that $\lambda \left( x^{\star }\right) =\lambda
\left( x\right) ^{\dagger },\qquad $ $\lambda ^{\bullet }\left( x^{\star
}\right) =\lambda _{\bullet }\left( x\right) ^{\dagger },\qquad $ $\lambda
_{\bullet }^{\bullet }\left( x^{\star }\right) =\lambda _{\bullet }^{\bullet
}\left( x\right) ^{\dagger }$
\end{theorem}

\proof
Let us represent the pre-Hilbert space $\mathfrak{D}$ as the $\mathfrak{D}_t$%
-span 
\begin{equation*}
\left\{ \sum_f\zeta ^f\otimes f^{\otimes }\left| \zeta ^f\in \mathfrak{D}%
_t,\quad f^{\bullet }\in \mathcal{K}\otimes L^2[t,\infty )\right. \right\}
\end{equation*}
of coherent vectors $f^{\otimes }\left( \tau \right) =\bigotimes_{t\in \tau
}f^{\bullet }\left( t\right) $, given for each finite subset $\tau =\left\{
t_1,...,t_n\right\} \subseteq [t,\infty )$ by tensor products $f^{\bullet
}\left( t_1\right) \otimes ...\otimes f^{\bullet }\left( t_n\right) $, with $%
\zeta ^f=0$ for almost all $f^{\bullet }$. Then the positivity (\ref{2.1})
of the quantum stochastic adapted maps $\phi _s,s>t$ into the $\mathfrak{D}$%
-forms $\left\langle h\right| \left. \phi _t\left( x\right) h\right\rangle ,$
for $h\in \mathfrak{D}$ can be obviously written as 
\begin{equation}
\sum_{x,y}\sum_{f,g}\left\langle \zeta _y^g\right| \left. \phi _s\left(
g^{\bullet },y^{\star }x,f^{\bullet }\right) \zeta _x^f\right\rangle \geq
0,\qquad  \label{2.7}
\end{equation}
the positive definiteness of the $\mathfrak{D}_t$-forms 
\begin{equation*}
\left\langle \eta \right| \left. \phi _s\left( g^{\bullet },x,f^{\bullet
}\right) \zeta \right\rangle =\left\langle \eta \otimes g^{\otimes }\right|
\left. \phi _s\left( x\right) \zeta \otimes f^{\otimes }\right\rangle
e^{-\int_s^\infty g^{\bullet }\left( r\right) ^{*}f^{\bullet }\left(
r\right) \mathrm{d}r},
\end{equation*}
with $\zeta _x^f\neq 0$ for a finite sequence of $x_k\in B,$ and for a
finite sequence of $f_k^{\bullet }\in \mathcal{K}\otimes L^2[t,\infty )$. If
the $\mathfrak{D}$-form $\phi _t\left( x\right) $ has the stochastic
differential (\ref{2.2}), the $\mathfrak{D}_t$-form $\phi _s\left(
g^{\bullet },x,f^{\bullet }\right) $ has the derivative 
\begin{eqnarray}
\frac{\mathrm{d}}{\mathrm{d}t}\phi _t\left( g^{\bullet },x,f^{\bullet
}\right) &=&g^{\bullet }\left( t\right) ^{*}f^{\bullet }\left( t\right) \phi
_t\left( x\right) +g^{\bullet }\left( t\right) ^{*}\alpha _{\bullet
}^{\bullet }\left( t,x\right) f^{\bullet }\left( t\right)  \label{2.8} \\
&&+g^{*}\left( t\right) ^{*}\alpha _{+}^{\bullet }\left( t,x\right) +\alpha
_{\bullet }^{-}\left( t,x\right) f^{\bullet }\left( t\right) +\alpha
_{+}^{-}\left( t,x\right)  \notag
\end{eqnarray}
at $s=t$. The positive definiteness, (\ref{2.7}), ensures the conditional
positivity 
\begin{equation}
\sum_{x,y}\sum_{f,g}\langle \zeta _y^g|\phi _t\left( y^{\star }x\right)
\zeta _x^f\rangle =0\Rightarrow \sum_{x,y}\sum_{f,g}\left\langle \zeta
_y^g\right| \left. \lambda _s^t\left( g^{\bullet },y^{\star }x,f^{\bullet
}\right) \zeta _x^f\right\rangle \geq 0  \label{2.9}
\end{equation}
of the form $\lambda _s^t\left( g^{\bullet },x,f^{\bullet }\right) =\frac
1{s-t}\left( \phi _s\left( g^{\bullet },x,f^{\bullet }\right) -\phi _t\left(
x\right) \right) $ for each $s>t$ and of the limit (\ref{2.8}) at $s\searrow
t$, given by the quadratic $\mathcal{K}$- form 
\begin{equation}
\lambda \left( g^{\bullet },x,f^{\bullet }\right) =b_{\bullet }\lambda
_{\bullet }^{\bullet }\left( x\right) a^{\bullet }+b_{\bullet }\lambda
^{\bullet }\left( x\right) +\lambda _{\bullet }\left( x\right) a^{\bullet
}+\lambda \left( x\right) ,  \label{2.10}
\end{equation}
where $a^{\bullet }=f^{\bullet }\left( t\right) ,\quad b_{\bullet
}^{*}=g^{\bullet }\left( t\right) $, and the $\lambda $'s are defined in (%
\ref{2.5}). Hence the form 
\begin{equation*}
\sum_{x,y}\sum_{\mu ,\nu }\left\langle \eta _y^\mu \right| \lambda _\nu ^\mu
\left( y^{\star }x\right) \eta _x^\nu \rangle :=\sum_{x,y}\left\langle \eta
_y\right| \lambda \left( y^{\star }x\right) \left| \eta _x\right\rangle
\end{equation*}
\begin{equation*}
+\sum_{x,y}\left( \left\langle \eta _y\right| \left. \lambda _{\bullet
}\left( y^{\star }x\right) \eta _x^{\bullet }\right\rangle +\left\langle
\eta _y^{\bullet }\right| \left. \lambda ^{\bullet }\left( y^{\star
}x\right) \eta _x\right\rangle +\left\langle \eta _y^{\bullet }\right|
\left. \lambda _{\bullet }^{\bullet }\left( y^{\star }x\right) \eta
_x^{\bullet }\right\rangle \right)
\end{equation*}
with $\eta =\sum_f\zeta ^f,\quad \eta ^{\bullet }=\sum_f\zeta ^f\otimes
f^{\bullet }\left( t\right) $ is positive if $\sum_{x,y}\langle \eta _y|\phi
_t\left( y^{\star }x\right) \eta _x\rangle =0$. The components $\eta \in 
\mathfrak{D}_t$ and $\eta ^{\bullet }\in \mathfrak{D}_t^{\bullet }$ are
independent because for any $\boldsymbol{\eta }\in \mathfrak{D}_t\oplus 
\mathfrak{D}_t^{\bullet }$ there exists such a function $a^{\bullet }\mapsto
\zeta ^a$ on $\mathcal{K}$ with a finite support, that $\sum_a\zeta ^a=\eta
,\quad \sum_a\zeta ^a\otimes a^{\bullet }=\eta ^{\bullet },$ namely, $\zeta
^a=0$ for all $a^{\bullet }\in \mathcal{K}$ except $a^{\bullet
}=a_n^{\bullet }$, the $n$-th $\mathcal{K}$-coefficient of the span $\eta
^{\bullet }=\sum_n\eta _n\otimes a_n^{\bullet }$, for which $\zeta ^a=\eta
_n $, and $a^{\bullet }=0$, for which $\zeta ^a=\eta -\sum_n\eta _n.$ This
proves the complete positivity of the matrix form $\boldsymbol{\lambda }$,
with respect to the matrix representation $\boldsymbol{\iota }$ defined in (%
\ref{2.5}) on the column-vectors $\boldsymbol{\eta }$. 
\endproof%

2. The conditional positivity of the germ-matrix $\boldsymbol{\lambda }$ at $%
t=0$ with respect to the representation $i:B\rightarrow \mathcal{A}\left( 
\mathcal{D}\right) $ embedded into the matrix form as in (\ref{2.5}),
obviously implies the positivity of the dissipation form 
\begin{equation}
\sum_{x,y}\left\langle \boldsymbol{\eta }_y\right| \boldsymbol{\Delta }%
\left( x,y\right) \left. \boldsymbol{\eta }_x\right\rangle
:=\sum_{x,y}\sum_{\mu ,\nu }\left\langle \eta _y^\mu \right| \Delta _\nu
^\mu \left( y,x\right) \left. \eta _x^\nu \right\rangle ,\qquad  \label{3.1}
\end{equation}
with $\eta =h_0\in \mathcal{D}$, $\eta ^{\bullet }=h_0^{\bullet }\in 
\mathcal{D}^{\bullet }$, $\eta ^{-}=\eta =\eta ^{+}$, corresponding to
non-zero $\boldsymbol{\eta }_x\in \mathcal{D}\oplus \mathcal{D}^{\bullet }$.
Here $\boldsymbol{\Delta }=\left( \Delta _\nu ^\mu \right) _{\nu =+,\bullet
}^{\mu =-,\bullet }$ is the stochastic dissipator, given by the blocks 
\begin{eqnarray}
\Delta _{\bullet }^{\bullet }\left( y,x\right) &=&\alpha _{\bullet
}^{\bullet }\left( y^{\star }x\right) +i\left( y^{\star }x\right) \otimes 
\mathbf{1},  \label{3.2} \\
\Delta _{\bullet }^{-}\left( y,x\right) &=&\alpha _{\bullet }^{-}\left(
y^{\star }x\right) -i\left( y\right) ^{*}\alpha _{\bullet }^{-}\left(
x\right) =\Delta _{+}^{\bullet }\left( x,y\right) ^{*}  \notag \\
\Delta _{+}^{-}\left( y,x\right) &=&\alpha _{+}^{-}\left( y^{\star }x\right)
-i\left( y\right) ^{*}\alpha _{+}^{-}\left( x\right) -\alpha _{+}^{-}\left(
y^{\star }\right) i\left( x\right) +i\left( y\right) ^{*}Di\left( x\right) ,
\notag
\end{eqnarray}
where $D=\lambda \left( 1\right) $. This means that the map $\lambda
_{\bullet }^{\bullet }$ is positive definite, and the conditions of the next
theorem are fulfilled if the maps $\lambda ,$ $\lambda ^{\bullet },$ $%
\lambda _{\bullet }$ have the following form 
\begin{eqnarray}
\lambda ^{\bullet }\left( x\right) &=&\varphi ^{\bullet }\left( x\right)
-K_{\bullet }^{\dagger }i\left( x\right) ,\qquad \text{ }\lambda _{\bullet
}\left( x\right) =\varphi _{\bullet }\left( x\right) -i\left( x\right)
K_{\bullet }\qquad  \label{3.3} \\
\lambda \left( x\right) &=&\varphi \left( x\right) -K^{\dagger }i\left(
x\right) -i\left( x\right) K,\qquad \text{ }\lambda _{\bullet }^{\bullet
}\left( x\right) =\varphi _{\bullet }^{\bullet }\left( x\right)  \notag
\end{eqnarray}
where $\boldsymbol{\varphi }=\left( \varphi _\nu ^\mu \right) _{\nu
=+,\bullet }^{\mu =-,\bullet }$ is a positive definite map from $B$ into the
sesquilinear forms on the pre-Hilbert space $\mathcal{D}\oplus \mathcal{D}%
^{\bullet }.$

\begin{theorem}
The following are equivalent:

\begin{enumerate}
\item[(i)] The dissipation form (\ref{3.1}), defined by the $\flat $-map $%
\alpha $ with $\alpha _{+}^{-}\left( 1\right) =D$, is positive definite: $%
\sum_{x,y}\left\langle \boldsymbol{\eta }_y\right| \boldsymbol{\Delta }%
\left( y,x\right) \left. \boldsymbol{\eta }_x\right\rangle \geq 0.$

\item[(ii)] There exists a pre-Hilbert space $\mathcal{D}^{\circ }$, a
unital $*$- representation $j$ of $B$ in $\mathcal{A}\left( \mathcal{D}%
^{\circ }\right) $, 
\begin{equation}
j\left( y^{\star }x\right) =j\left( y\right) ^{*}j\left( x\right) ,\quad
j\left( 1\right) =I,  \label{3.6}
\end{equation}
a $\left( i,j\right) $-derivation of $B$, 
\begin{equation}
k\left( y^{\star }x\right) =j\left( y\right) ^{*}k\left( x\right) +k\left(
y^{\star }\right) i\left( x\right) ,  \label{3.7}
\end{equation}
having values in the operators $\mathcal{D}\rightarrow \mathcal{D}^{\circ }$
, the adjoint map $k^{*}\left( x\right) =k\left( x^{\star }\right) ^{\dagger
}$, 
\begin{equation*}
k^{*}\left( y^{\star }x\right) =i\left( y\right) ^{*}k^{*}\left( x\right)
+k^{*}\left( y^{\star }\right) j\left( x\right)
\end{equation*}
into the operators $\mathcal{D}^{\circ }\rightarrow \mathcal{D}^{\prime }$,
and a map $l:B\rightarrow \mathcal{Q}\left( \mathcal{D}\right) $, having the
coboundary property 
\begin{equation}
l\left( y^{\star }x\right) =i\left( y\right) ^{*}l\left( x\right) +l\left(
y^{\star }\right) i\left( x\right) +k^{*}\left( y^{\star }\right) k\left(
x\right) ,  \label{3.8}
\end{equation}
such that $l\left( x\right) +Di\left( x\right) =\lambda \left( x\right)
=l^{*}\left( x\right) +i\left( x\right) ^{\prime }D$, where $l^{*}\left(
x\right) =l\left( x^{\star }\right) ^{\dagger }$, 
\begin{equation*}
L_{\circ }^{\bullet }k\left( x\right) +L_{+}^{\bullet }i\left( x\right)
=\lambda ^{\bullet }\left( x\right) ,\quad \lambda _{\bullet }\left(
x\right) =k^{*}\left( x\right) L_{\bullet }^{\circ }+i\left( x\right)
^{\prime }L_{\bullet }^{-},
\end{equation*}
and $\lambda _{\bullet }^{\bullet }\left( x\right) =L_{\circ }^{\bullet
}j\left( x\right) L_{\bullet }^{\circ }$ for some operators $L_{+}^{\bullet
}:\mathcal{D}\rightarrow \mathcal{D}^{\bullet }$, $L_{\circ }^{\bullet }:%
\mathcal{D}^{\circ }\rightarrow \mathcal{D}^{\bullet }$, with the adjoints $%
L_{\bullet }^{-}=L_{+}^{\bullet \dagger }$, $L_{\bullet }^{\circ }=L_{\circ
}^{\bullet \dagger }$.

\item[(iii)] There exists a pseudo-Hilbert space, $\mathcal{E}$, a unital $%
\flat $-representation $\jmath :B\mathcal{\rightarrow A}\left( \mathcal{E}%
\right) $, and a linear operator $L:\mathcal{D}\oplus \mathcal{D}^{\bullet
}\rightarrow \mathcal{E}$, where $\mathcal{D}^{\bullet }=\mathcal{K}\otimes 
\mathcal{D}$ such that 
\begin{equation}
\boldsymbol{\lambda }\left( x\right) =L^{\flat }\jmath \left( x\right)
L,\qquad \forall x\in B.  \label{3.10}
\end{equation}

\item[(iv)] The germ-mapping $\boldsymbol{\lambda }=\boldsymbol{i}+%
\boldsymbol{\alpha }$ is conditionally positive definite with respect to the
matrix representation $\boldsymbol{\iota }$ in (\ref{2.5}), where $\iota =i$.
\end{enumerate}
\end{theorem}

\proof
The proof of the implication (i)$\Rightarrow $(ii), generalizing the
Evans-Lewis Theorem, is similar to the proof of the dilation theorem in \cite%
{Bcs}.The proof of the implication (ii)$\Rightarrow $(iii) can be also
obtained as in \cite{Bpe} by the explicit construction of $\mathcal{E}$ as $%
\oplus _{\mu =-}^{+}\mathcal{D}^\mu $, where $\mathcal{D}^{+}=\mathcal{D}$, $%
\mathcal{D}^{-}=\mathcal{D}^{\prime }$, with the indefinite metric tensor $%
\mathbf{G}=\left[ G_{\mu \nu }\right] $ given above for $\mu ,\nu =-,\circ
,+ $, and $D=\lambda \left( 1\right) ,$ of the unital $\flat $%
-representation $\mbox{\boldmath $\j $}=\left[ \jmath _\nu ^\mu \right]
_{\nu =-,\circ ,+}^{\mu =-,\circ ,+}$ of $B$ on $\mathcal{E}$ : 
\begin{equation*}
\mbox{\boldmath $\j $}\left( y^{\star }x\right) =\mbox{\boldmath $\j $}%
\left( y\right) ^{\flat }\mbox{\boldmath $\j $}\left( x\right) ,\quad %
\mbox{\boldmath $\j $}\left( 1\right) =\mathbf{I}
\end{equation*}
with $\mbox{\boldmath $\j $}\left( x\right) ^{\flat }=\mathbf{G}^{-1}%
\mbox{\boldmath $\j $}\left( x\right) ^{\dagger }\mathbf{G}$, given by the
components 
\begin{equation}
\jmath _{+}^{+}=i,\quad \jmath _{\circ }^{\circ }=j,\quad \jmath
_{-}^{-}=i^{\prime }\quad \jmath _{+}^{\circ }=k,\quad \jmath _{\circ
}^{-}=k^{*},\quad \jmath _{+}^{-}=l,  \label{3.9}
\end{equation}
where $i^{\prime }\left( x\right) =i\left( x\right) ^{\prime }$ and all
other $\jmath _\nu ^\mu =0$. The linear operator $L$ is given by $\mathbf{L=}%
\left[ L^\mu ,L_{\bullet }^\mu \right] $ with the components 
\begin{equation*}
L^{+}=I,\quad L^{\circ }=0,\quad L^{-}=0,\quad L_{\bullet }^{+}=0,\quad
L_{\bullet }^{\circ }=L_{\circ }^{\bullet \dagger },\quad L_{\bullet
}^{-}=L_{+}^{\bullet \dagger },
\end{equation*}
and $\mathbf{L}^{\flat }=\left( 
\begin{array}{ccc}
I & 0 & D \\ 
0 & L_{\circ }^{\bullet } & L_{+}^{\bullet }%
\end{array}
\right) =\mathbf{L}^{\dagger }\mathbf{G}$ such that $\mathbf{L}^{\flat }%
\mbox{\boldmath $\j $}\left( x\right) \mathbf{L}=\boldsymbol{\lambda }\left(
x\right) $. The implication (iii)$\Rightarrow $(iv) is a straight forward
consequence of this construction, and the implication (iv)$\Rightarrow $(i)
is similar to the non-stochastic case \cite{Lnd} . 
\endproof%

\end{document}